%% file: notes_nov-09.tex
\documentclass[11pt,oneside]{article}
\usepackage{amssymb,eucal}
\usepackage{amsfonts,amsmath,epsfig}

\usepackage[active]{srcltx}
\usepackage[sort&compress,square]{natbib}
\usepackage[colorlinks]{hyperref}

\input myamsdefs

\numberwithin{equation}{section}
\newtheorem{theorem}{Theorem}[section]
\newtheorem{prop}[theorem]{Proposition}
\newtheorem{lemma}[theorem]{Lemma}
\newtheorem{assumption}{Condition}

\def\qed{\hfill$\Box$\medbreak}

\newenvironment{proof}[1][Proof]
{
  \medbreak\noindent\textit{#1.\/}
}{
  \qed
}

\newenvironment{proof*}[1][Proof]
{
  \medbreak\noindent\textit{#1.\/}
}{
}

\long\def\comment#1{{}}

\def\r{\epsilon}

\def\normx#1{\|#1\|}
\def\angx#1#2{\langle{#1},\,{#2}\rangle}
\def\tp{\sp{\top}}

\def\mean{\text{\sf E}}

\def\prob{\text{\sf Pr}}
\def\sppt{\mathrm{spt}}

\def\argx#1{\mathop{\arg#1}}

\def\rc{\varrho}  

\def\pset#1#2{\mathcal{D}\Grp{#1,#2}}  

\def\parsp{{\Reals^p}}    
\def\obssp{{\Reals^n}}    
\def\est#1{\widehat{#1}}
\def\btw#1#2#3{{#2\le#1\le#3}}

\def\diam{\delta}

\def\sfd{\text{\sf d}}

\def\cD{\mathcal{D}}

\def\cN{\mathcal{N}}

\begin{document}
\begin{center}
  \large
  \bf
  On $\vf{\ell_1}$-regularized estimation for nonlinear models that
  have sparse underlying linear structures  \\[1.5ex]
  \rm\normalsize
  Zhiyi Chi\footnote{Address:  215 Glenbrook Road, U-4120, Storrs, CT
    06269, USA}
  \\
  Department of Statistics \\
  University of Connecticut \\[1.5ex]
  \today
\end{center}

\begin{abstract}
  In \cite{chi:09}, for nonlinear models with sparse underlying linear
  structures, we studied the error bounds of $\ell_0$-regularized
  estimation.  In this note, we show that $\ell_1$-regularized
  estimation in some important cases can achieve the same order of error
  bounds as those in \cite{chi:09}.

  \medbreak\noindent
  \textit{Keywords and phrases.\/} Regularization, sparsity, MLE,
  regression, variable selection, parameter estimation, nonlinearity,
  power series expansion, analytic, exponential.

  \medbreak\noindent
  \textit{AMS 2000 subject classification.\/} Primary 62G05;
  secondary 62J02.

  \medbreak\noindent
  \textit{Acknowledgement.\/}  Research partially supported by
  NSF Grant DMS-07-06048 and NIH Grant MH-68028.

\end{abstract}

\section{Introduction} \label{sec:intro}
The models we consider are of the form
\begin{align} \label{eq:nonlinear}
  y = f(X\tp \beta) + \r, \quad X\in\Reals^{n\times p}, \ \
  \beta\in\parsp,\ \ y,\ \r \in\obssp,
\end{align}
where $f: \Reals\to\Reals$ is a \emph{known\/} function, $X$ a fixed
design matrix, $y$ and $\r$ are vectors of observations and errors,
respectively.  In \eqref{eq:nonlinear} and henceforth, for
$x\in\obssp$, we denote $f(x) = (f(x_1), \ldots, f(x_n))\tp$.
The parameter $\beta$ is sparse in the sense that the number of
its nonzero coordinates is much smaller than its dimension
\cite{wasserman:09}.

For $a>0$ and $v\in\parsp$, denote by $\normx{v}_a$ the $\ell_a$-norm
of $v$.  The support of $v$ is defined to be $\sppt(v) :=\{i:
v_i\not=0\}$.  Denote by $|A|$ the cardinality of a set $A$.  The
$\ell_0$-norm of $v$ is $\normx{v}_0=|\sppt(v)|$.  By an
$\ell_a$-regularized estimator of $\beta$ we mean
\begin{align} \label{eq:general} 
  \est\beta = \argx\min_{v\in D} \Sbr{\ell(y, Xv) + c_r\normx{v}_a},
\end{align}
where $D\subset\parsp$ is a pre-selected search domain, $\ell(y, Xv)$
a loss function, and $c_r>0$ a tuning parameter.  We are interested in
the case where $a=1$.

For models \eqref{eq:nonlinear}, much has been learned about the case
where $p$ is fixed or much smaller than $n$ (cf.\ \cite{hristache:01,
  hristache:01b} and references therein).  The note is concerned with
the case where $p$ can be large, possibly much larger than $n$ and, at
the same time, $|\sppt(\beta)|$ is much smaller than $p$.  Under this
setting, the case where $f(x)=x$ has been a subject of great interest
recently (cf.\ \citep{tibshirani:96, zhang:09, candes:plan:09,
  zhao:yu:06, bunea:etal:07, candes:tao:07} and references therein).

The main purpose of the note is to establish general results on the
estimator \eqref{eq:general} similar to Proposition~2.1 in
\cite{chi:09}.  Once established, the results allow the steps in
\cite{chi:09} to be followed, often word by word, to get error bounds
for specific cases.  In \eqref{eq:general},  while the function being
maximized only involves $\normx{v}_1$, the search domain $D$ may be
constrained in terms of $\normx{v}_0$ as well as certain weighted
$\ell_1$-norm of $v$.  As a result, we get two types of estimators,
one being regularized by $\normx{v}_0$ and (weighted) $\ell_1$-norms
of $v$, the other only by $\ell_1$-norms of $v$.  Error bounds for
both types of estimators will be derived.  The former type of
estimators can attain the same order of precision as their
$\ell_0$-regularized counterparts studied in \cite{chi:09}.  In
contrast, although the latter type of estimators are computationally
more amenable, in some cases they seem unable to attain the same order
of precision, at least with the techniques employed here.

To reduce repetition, we will omit most of results that can be
established directly following \cite{chi:09} and instead focus on
those that require new ideas.

\section{Main results} \label{sec:main}
The row vectors and column vectors of $X$ will be denoted by $X_1\tp$,
\ldots, $X_n\tp$ and $\eno V p$, respectively.  We shall assume that
$V_j\not=0$.  For $g=(g_1, \ldots, g_n)$ and $x\in \obssp$, where each
$g_i: \Reals\to\Reals$ is a function, denote $g(x) = (g_1(x_1),
\ldots, g_n(x_n))\tp$.

As in \cite{chi:09}, to bound the error of the $\ell_1$-regularized
estimator in \eqref{eq:general}, our first step is to show that
$\est\beta$ belongs 
to a set of $v$ that satisfy the following inequality,
\begin{align} \label{eq:ineq}
  G(\psi(X v) - \psi(X\beta))
  \le 2 |\angx{\r}{\varphi(Xv) - \varphi(X\beta)}|
  -
  c_r(\normx{v}_1-\normx{\beta}_1),
\end{align}
where $G: \obssp\to\Reals$ is a function, $\psi = (\eno\psi n)$, 
$\varphi = (\eno \varphi n)$, with $\psi_i$ and $\varphi_i$ being
functions from $\Reals$ to $\Reals$.  In many cases, it is not very
hard to get \eqref{eq:ineq} for maximum likelihood estimators (MLE) or
least square estimators (LSE).  We will illustrate this later.  Our
focus next is to use \eqref{eq:ineq} to derive two error bounds for
$\est\beta$.

\subsection{Conditions and general error bounds}
\renewcommand{\theassumption}{H\arabic{assumption}}
For both error bounds, we need the following condition.
\begin{assumption}\rm\label{cond:A1}
  Given $q\in (0,1)$, there is $c_1=c_1(X, \beta, \varphi,q)>0$, such
  that
  \begin{align*}
    \prob\Cbr{
      |\angx{\r}{\varphi(Xv) - \varphi(X\beta)}|
      \le c_1\sqrt{n} \normx{v-\beta}_1,
      \text{ all } v \in D
    }\ge 1-c_0 q,
  \end{align*}
  where $c_0>0$ is an arbitrarily pre-selected constant, such as 1 or 2.
\end{assumption}

The same condition was used in \cite{chi:09}, but with $c_0=2$.
As remarked in \cite{chi:09}, $c_0$ is purely for notational
ease when Condition \ref{cond:A1} is verified for specific
cases.  To get the first bound, we also need another condition used in
\cite{chi:09}.
\begin{assumption}\rm \label{cond:A2}
  There is $c_2=c_2(X, \beta, \psi)>0$, such that for \emph{all\/}
  $v\in D$,
  \begin{align*}
    G(\psi(X v) - \psi(X\beta)) \ge c_2 n \normx{v-\beta}_2^2.
  \end{align*}
\end{assumption}

We now can state the first error bound for $\est\beta$.
\begin{prop} \label{prop:basic}
  Suppose \emph{Conditions \ref{cond:A1} and \ref{cond:A2}\/} are
  satisfied.  If $\est\beta\in D$ is a random variable that always
  satisfies the inequality \eqref{eq:ineq} with $c_r = 2 c_1
  \sqrt{n}$, then, letting $\kappa_r = 4 c_1/c_2$,
  $\prob\{\normx{\est\beta-\beta}_2 \le \kappa_r
  \sqrt{|\sppt(\beta)|/n}\} \ge 1-c_0 q$.
\end{prop}

To get the second bound, we replace Condition \ref{cond:A2} with the
next one.
\begin{assumption} \rm \label{cond:A3}
  There is $c_3 = c_3(X,\beta,\psi)>0$, such that for \emph{all\/}
  $z\in \{X v: v\in D\}$,
  \begin{align*}
    G(\psi(z) - \psi(X\beta)) \ge c_3 \normx{z-X \beta}_2^2.
  \end{align*}
\end{assumption}

We also need some conditions on the second moments of the column
vectors of $X$.  Such conditions are sometimes referred to as
coherence property \citep{bunea:etal:07, candes:plan:09}.  Let
\begin{align*}
  \mu_X = \max_{1\le i<j\le p} \frac{|V_i\tp V_j|}{\normx{V_i}_2
    \normx{V_j}_2}, \quad
  a_X = \min_{\btw i 1 p} \frac{\normx{V_i}_2^2}{n}, \quad
  b_X = \max_{\btw i 1 p} \frac{\normx{V_i}_2^2}{n}.
\end{align*}

\begin{prop} \label{prop:basic2}
  Suppose \emph{Conditions \ref{cond:A1}\/} and \emph{\ref{cond:A3}\/}
  are satisfied and $a_X + b_X\mu_X > 6 b_X |\sppt(\beta)| \mu_X$.
  Fix $\tau>0$ such that
  \begin{align} \label{eq:t}
    a_X + b_X\mu_X > 2 b_X(3 + 4\tau)|\sppt(\beta)|\mu_X,
  \end{align}
  and let $c_r = 2 (1+1/\tau) c_1\sqrt{n}$,
  \begin{align*}
    \kappa_r = \frac{3(2+1/\tau) \sqrt{2+(1+2\tau)^2}}
    {a_X+b_X\mu_X}  \times  \frac{c_1}{c_3}.
  \end{align*}
  If $\est\beta\in D$ is a random variable that always satisfies the
  inequality \eqref{eq:ineq} with the above $c_r$ as the tuning
  parameter, then $\prob\{\normx{\est\beta-\beta}_2 \le \kappa_r
  \sqrt{|\sppt(\beta)|/n}\} \ge 1-c_0 q$.
\end{prop}

Since $a_X\le b_X$, \eqref{eq:t} sets an upper bound on $\mu_X$.  To
get a moderate value of $\kappa_r$ in Proposition \ref{prop:basic2},
$\tau$ has to be moderate.  If, say, $\tau=1$, then by \eqref{eq:t},
$a_X/b_X > (14|\sppt(\beta)|-1) \mu_X$, which further limits the
magnitude of $\mu_X$.  Under certain conditions, one can get
$\mu_X=O(\sqrt{n^{-1}\ln p})$ \cite{candes:plan:09, chi:09}, which is
small for large $n$, even when $p$ is much larger than $n$, for
example, $p=n^\alpha$ with some $\alpha>1$.

We next make some comments on conditions used in specific cases to
establish Conditions \ref{cond:A1} -- \ref{cond:A3}.  To establish
Condition \ref{cond:A1}, the following tail condition on the errors
$\r_i$ is useful: there are $\sigma>0$ and $c_\r\ge 1$,  such that
\begin{align} \label{eq:tail}
  \prob\{|a\tp \r|>t\normx{a}_2\} \le c_\r e^{-t^2/(2\sigma^2)},
  \quad
  \text{all}\ \ t\ge 0, \ a\in\obssp.
\end{align}
As remarked in \cite{chi:09}, typically $c_\r$ can be set at 2.
At the end of the note, we will see that in some cases $c_\r$
has to be set at other values.

To establish Condition \ref{cond:A2} or \ref{cond:A3}, we usually need
to put some restrictions on the search domain $D$ in
\eqref{eq:general}.  To establish Condition \ref{cond:A3}, which is
the less restrictive of the two, we typically choose
\begin{align}
  D\subseteq\cD(I)
  = T^{-1}(I^n) = 
  \{v\in\parsp: X_i\tp v\in I, \ \btw i 1 n\},  \label{eq:domain}
\end{align}
where $T$ is the mapping $v\to Xv$ and $I$ is an interval in $\Reals$.
In general, we need not put restrictions on $|\sppt(v)|$.
\comment{Note that all the constraints in \eqref{eq:domain} are linear
  and $\cD(I)$ is convex.}  On the other hand, to establish Condition
\ref{cond:A2}, we typically start with verifying Condition
\ref{cond:A3}, and then proceed to get $\normx{X(v-\beta)}_2\ge c
\normx{v-\beta}_2$ for some constant $c>0$.  To do this, we need to
put restrictions on $|\sppt(v)|$, typically by
requiring \begin{align*}   
  D\subseteq\pset I h
  =\cD(I) \cap \Cbr{u\in \parsp: |\sppt(u)|\le h},
\end{align*}
with $h\ge 1$ being bounded in terms of $\mu_X$ (cf.\ \cite{chi:09}).
Thus, though not directly used in Proposition \ref{prop:basic}, 
coherence property of $X$ is needed in specific applications of the
Proposition.

\subsection{Proofs} \label{ssec:proof}

For $v\in\parsp$ and $S\subset\{1,\ldots,p\}$, denote $v_S=(\eno x
p)^T$ with $x_i= v_i \cf{i\in S}$.  Let $d=v-\beta$.  Then for any
$S\supset\sppt(\beta)$, we have $v = \beta+d_S +v_{S^c}$ and 
\begin{align} \label{eq:norm-decom}
  \normx{v}_1
  = \normx{\beta+d_S}_1 + \normx{v_{S^c}}_1, \quad
  \normx{d}_a^a
  = \normx{d_S}_a^a + \normx{v_{S^c}}_a^a, \quad
  \text{for any}\ a>0.
\end{align}

\begin{proof}[Proof of Proposition \ref{prop:basic}]
  Let $d=\est\beta-\beta$.  Because $\est\beta$ always
  satisfies \eqref{eq:ineq}, by Conditions \ref{cond:A1} and
  \ref{cond:A2}, with probability at least $1-c_0 q$,
  \begin{align*}
    c_2 n \normx{d}_2^2
    &\le 2 c_1\sqrt{n}\normx{d}_1
    - c_r(\normx{\est\beta}_1 - \normx{\beta}_1) \\
    &= 2 c_1\sqrt{n} (\normx{d}_1 + \normx{\beta}_1 -
    \normx{\est\beta}_1).
  \end{align*}
  Let $S=\sppt(\beta)$.  Apply \eqref{eq:norm-decom} to the right hand
  side of the above inequality to get
  \begin{align*}
    c_2 n \normx{d}_2^2
    &
    \le 2c_1\sqrt{n} ( \normx{d_S}_1 + 
    \normx{\est\beta_{S^c}}_1 + \normx{\beta}_1 - 
    \normx{\beta+ d_S}_1 - \normx{\est\beta_{S^c}}_1)\\
    &
    = 2c_1\sqrt{n} ( \normx{d_S}_1 + 
    \normx{\beta}_1 - \normx{\beta+ d_S}_1).
  \end{align*}
  Then by Minkowski inequality and Cauchy-Schwartz inequality, 
  \begin{align*}
    \normx{d}_2^2
    \le 4(c_1/c_2)\normx{d_S}_1/\sqrt{n}
    \le \kappa_r \sqrt{|S|/n}\,\normx{d_S}_2.
  \end{align*}
  Because $\normx{d}_2^2 = \normx{d_S}_2^2 +
  \normx{\est\beta_{S^c}}_2^2$ by \eqref{eq:norm-decom}, the above
  inequalities imply
  $$
  \normx{d}_2^2 \le
  M:=\sup\{x^2 + y^2: \text{$x\ge 0$ and $y\ge 0$ satisfy}\ x^2 + y^2
  \le \kappa_r \sqrt{|S|/n} x\}.
  $$
  
  To find $M$, first, in order that $x^2 + y^2 \le \kappa_r
  \sqrt{|S|/n} x$, there must be $\kappa_r^2 |S|/n \ge
  4y^2$.  Given $y\ge 0$ satisfying the condition, the maximum
  possible $x$ is 
  $$
  x_0(y) = (1/2)[\kappa_r \sqrt{|S|/n} + \sqrt{\kappa_r^2 |S|/n -
      4y^2}\,].
  $$
  It is seen that
  $$
  x_0^2(y) + y^2 = \frac{\kappa_r^2 |S|/n + \kappa_r\sqrt{|S|/n}
    \sqrt{ \kappa_r^2 |S|/n - 4y^2}}{2}
  \le \kappa_r^2 |S|/n.
  $$
  Therefore, $M=\kappa_r^2|S|/n$, where the maximum is obtained if and
  only if $x=\kappa_r \sqrt{|S|/n}$ and $y=0$.  This yields 
  $\normx{d}_2\le\sqrt{M}=\kappa_r \sqrt{|S|/n}$, as desired.
\end{proof}

\begin{proof}[Proof of Proposition \ref{prop:basic2}]
  It suffices to show that
  \begin{align} \label{eq:precision}
    \prob\Cbr{\normx{v-\beta}_2 \le
      \kappa_r\sqrt{|\sppt(\beta)|/n}
      \text{ for \emph{all\/} $v\in D$ satisfying \eqref{eq:ineq}}
    } \ge 1-c_0 q.
  \end{align}

  By Conditions \ref{cond:A1} and \ref{cond:A3}, with probability at
  least $1-c_0 q$, the inequality 
  \begin{align*}
    c_3 \normx{X (v-\beta)}_2^2
    \le 2 c_1\sqrt{n}\normx{v-\beta}_1 - 2
    c_1(1+1/\tau)\sqrt{n}(\normx{v}_1 - \normx{\beta}_1)
  \end{align*}
  holds for \emph{all\/} $v\in D$ satisfying \eqref{eq:ineq}.  Fix one
  such $v$ and an arbitrary $S\supset\sppt(\beta)$.  Let $d=v-\beta$.
  By \eqref{eq:norm-decom},
  \begin{align*}
    c_3 \normx{X d}_2^2
    &
    \le 2 c_1\sqrt{n}(\normx{d_S}_1+\normx{v_{S^c}}_1) \\
    &\hspace{2cm}
    - 2 c_1(1+1/\tau)\sqrt{n}(\normx{\beta+d_S}_1 +
    \normx{v_{S^c}}_1- \normx{\beta}) \\
    &
    = 2 c_1\sqrt{n}\normx{d_S}_1
    - 2 c_1(1+1/\tau)\sqrt{n}(\normx{\beta+d_S}_1 - \normx{\beta})
    - 2(c_1/\tau)\sqrt{n}\normx{v_{S^c}}_1.
  \end{align*}

  For ease of notation, denote $\tilde c_1 = c_1/c_3$ for now.  By
  Minkowski inequality, $\normx{\beta+d_S}_1-\normx{\beta}_1\ge
  -\normx{d_S}_1$, and so
  \begin{align} \label{eq:norm-1}
    \normx{X d}_2^2 
    \le 2 \tilde c_1(2+1/\tau) \sqrt{n}\normx{d_S}_1 -
    (2\tilde c_1/\tau) \sqrt{n}\normx{v_{S^c}}_1.
  \end{align}
  
  First of all, since the left hand side of \eqref{eq:norm-1} is
  nonnegative, it follows that
  \begin{align} \label{eq:v-1}
    \normx{v_{S^c}}_1 \le (1+2\tau) \normx{d_S}_1.
  \end{align}
  On the other hand, by $Xd=Xd_S + Xv_{S^c}$,
  \begin{align*}
    \normx{X d}_2^2
    &
    = \normx{X d_S}_2^2 + \normx{X v_{S^c}}_2^2 +
    2 \Ang{X d_S, X v_{S^c}}
    \ge
    \normx{X d_S}_2^2 - 2\Abs{\Ang{X d_S, X v_{S^c}}}.
  \end{align*}
  We next derive a lower bound of $\normx{X d}_2^2$.  First,
  by $X d_S = \sum_{i\in S} d_i V_i$, 
  \begin{align*}
    \normx{X d_S}_2^2
    &
    = \sum_{i\in S} d_i^2 \normx{V_i}_2^2
    + \sum_{i,j\in S,\, i\not=j} d_i d_j (V_i\tp V_j) \\
    &
    \ge \sum_{i\in S} d_i^2 \normx{V_i}_2^2
    - \sum_{i,j\in S,\,i\not=j} |d_i d_j| |V_i\tp V_j|.
  \end{align*}
  Because $\normx{V_i}_2^2\ge a_X$ and for $i\not=j$, $|V_i\tp V_j|\le
  \mu_X \normx{V_i}_2 \normx{V_j}_2 \le b_X \mu_X n$, we get
  \begin{align*}
    &\hspace{-1cm}
    \normx{ X d_S}_2^2 \ge a_X n \sum_{i\in S} d_i^2 - 
    b_X\mu_X n \sum_{i,j\in S,\, i\not=j}
    |d_i d_j| \\
    &
    = (a_X +b_X\mu_X) n \normx{d_S}_2^2 -
    b_X\mu_X n \normx{d_S}_1^2.
  \end{align*}
  Second, by $X v_{S^c} = \sum_{j\not\in S} v_j V_j$,
  \begin{align*}
    \Abs{\Ang{X d_S, X v_{S^c}}}
    &
    =\Abs{\sum_{i\in S,\ j\not\in S} d_i v_j V_i\tp V_j }
    \le \sum_{i\in S, \ j\not\in S} |d_i v_j| |V_i\tp V_j|\\
    &
    \le
    b_X\mu_X n \sum_{i\in S, \ j\not\in S} |d_i v_j| 
    = b_X\mu_X n \normx{d_S}_1 \normx{v_{S^c}}_1.
  \end{align*}
  Therefore, putting the above inequalities together,
  \begin{align} \label{eq:norm-2}
    &
    \normx{X d}_2^2
    \ge 
    (a_X +b_X\mu_X) n \normx{d_S}_2^2 -
    b_X\mu_X n \normx{d_S}_1^2 - 
    2 b_X\mu_X n \normx{d_S}_1 \normx{v_{S^c}}_1.
  \end{align}
  Combining \eqref{eq:norm-1} and \eqref{eq:norm-2}, and then grouping
  the terms, we get
  \comment{
    \begin{multline*}
      (a_X +b_X\mu_X) n \normx{d_S}_2^2 -
      b_X\mu_X n \normx{d_S}_1^2 - 
      2 b_X\mu_X n \normx{d_S}_1 \normx{v_{S^c}}_1 \\
      \le 2 \tilde c_1(2+1/\tau) \sqrt{n}\normx{d_S}_1 -
      (2\tilde c_1/\tau) \sqrt{n}\normx{v_{S^c}}_1,
    \end{multline*}
    and hence
  }
  \begin{align}
    (a_X +b_X\mu_X) n \normx{d_S}_2^2
    &
    \le
    b_X\mu_X n \normx{d_S}_1^2 +
    2 \tilde c_1(2+1/\tau) \sqrt{n}\normx{d_S}_1 \nonumber\\
    &
    \qquad+
    2 \Cbr{b_X\mu_X \sqrt{n} \normx{d_S}_1 - \tilde c_1/\tau}
    \sqrt{n} \normx{v_{S^c}}_1.
    \label{eq:norm-3}
  \end{align}

  So far, other than the requirement that $S\supset \sppt(\beta)$,
  the choice of $S$ is arbitrary.  To continue, we need the next result
  that puts more constraints on $S$.
  \begin{lemma} \label{lemma:bound-1}
    Suppose $S\supset\sppt(\beta)$ such that $a_X+b_X\mu_X > b_X
    \mu_X(3+4\tau)|S|$.  Then $b_X\mu_X \sqrt{n}
    \normx{d_S}_1 < \tilde c_1/\tau$.
  \end{lemma}

  Assume the lemma is true for now.  Let $S\supset\sppt(\beta)$ such
  that $a_X+b_X \mu_X >b_X\mu_X(3+4\tau)|S|$.  Later we will see that
  such $S$ indeed exists and make specific choices for it.  By
  \eqref{eq:norm-3}, Lemma \ref{lemma:bound-1}, and
  Cauchy-Schwartz inequality,
  \begin{align*}
    (a_X +b_X\mu_X) n \normx{d_S}_2^2 
    &\le b_X\mu_X n \normx{d_S}_1^2
    + 2 \tilde c_1(2+1/\tau) \sqrt{n}\normx{d_S}_1 \\
    &
    \le
    b_X\mu_X n |S|\normx{d_S}_2^2
    + 2 \tilde c_1(2+1/\tau) \sqrt{n |S|}\normx{d_S}_2 \\
    &
    \le
    (a_X+b_X\mu_X) n \normx{d_S}_2^2/3 + 2 \tilde c_1(2+1/\tau)
    \sqrt{n|S|}\normx{d_S}_2,
  \end{align*}
  where the last inequality is due to $a_X + b_X \mu_X > 3 b_X \mu_X
  |S|$.  Thus
  \begin{align} \label{eq:v-2}
    \normx{d_S}_2 \le
    \frac{3\tilde c_1(2+1/\tau)\sqrt{|S|}}{(a_X+b_X\mu_X)\sqrt{n}}.
  \end{align}
  
  Let $S_1$ be the union of $\sppt(\beta)$ and the set of $i\not\in
  \sppt(\beta)$ with the $|\sppt(\beta)|$ largest $d_i$ outside
  $\sppt(\beta)$.  By Lemma 3.1 of \cite{candes:tao:07},
  \begin{align} \label{eq:v-3}
    \normx{d}_2^2
    \le
    \normx{d_{S_1}}_2^2 +
    \frac{\normx{d_{\sppt(\beta)^c}}_1^2}{|\sppt(\beta)|}.
  \end{align}
  Since $d_{\sppt(\beta)^c}=v_{\sppt(\beta)^c}$, by \eqref{eq:v-1} and
  Cauchy-Schwartz inequality,
  $$
  \normx{d_{\sppt(\beta)^c}}_1
  \le (1+2\tau)\normx{d_{\sppt(\beta)}}_1
  \le (1+2\tau)\sqrt{|\sppt(\beta)|} \normx{d_{\sppt(\beta)}}_2,
  $$
  which together with \eqref{eq:v-3} yields
  \begin{align} \label{eq:v-4}
    \normx{d}_2^2 \le \normx{d_{S_1}}_2^2 +
    (1+2\tau)^2 \normx{d_{\sppt(\beta)}}_2^2.
  \end{align}

  Note $|S_1|=2|\sppt(\beta)|$.  By the assumption in \eqref{eq:t}
  and Lemma \ref{lemma:bound-1}, it is seen that \eqref{eq:v-2} holds
  for $S=S_1$ and for $S=\sppt(\beta)$.  Combine this with
  \eqref{eq:v-4} to get
  \begin{align*}
    \normx{d}_2 \le 
    \frac{3c_1(2+1/\tau)
      \sqrt{[2+(1+2\tau)^2]\,|\sppt(\beta)|}}
    {c_3(a_X+b_X\mu_X)\sqrt{n}}\,,
  \end{align*}
  where we have recovered $\tilde c_1=c_1/c_3$.  The proof of
  \eqref{eq:precision} is then complete. 
\end{proof}

\begin{proof}[Proof of Lemma \ref{lemma:bound-1}]  
  Assume the opposite were true, i.e.\ $b_X\mu_X \sqrt{n}
  \normx{d_S}_1 \ge \tilde c_1/\tau$.  Then 
  clearly $d_S\not=0$.   By \eqref{eq:v-1}, the right hand side of
  \eqref{eq:norm-3} is no greater than
  \begin{multline*}
    2 \tilde c_1(2+1/\tau) \sqrt{n}\normx{d_S}_1 +
    2 \Cbr{b_X\mu_X \sqrt{n} \normx{d_S}_1 - \tilde c_1/\tau}
    \sqrt{n} (1+2\tau)\normx{d_S}_1 \\
    = 2b_X\mu_X n (1+2\tau)\normx{d_S}_1^2,
  \end{multline*}
  so \eqref{eq:norm-3} together with Cauchy-Schwartz inequality yields
  $(a_X+b_X\mu_X) \normx{d_S}_2^2 \le b_X\mu_X (3+4\tau)
  \normx{d_S}_1^2 \le b_X \mu_X (3+4\tau) |S| \normx{d_S}_2^2$.  Since
  $d_S\not=0$, then $a_X+b_X\mu_X \le b_X \mu_X (3+4\tau) |S|$, which
  contradicts the assumption.
\end{proof}

\section{MLE for exponential linear models and LSE for analytic
  models} \label{sec:example}

In \cite{chi:09}, by choosing suitable search domain $D$, we derived
error bounds for the $\ell_0$-regularized MLE and LSE for exponential
linear models and analytic models, respectively.  Under the
conditions in Proposition \ref{prop:basic}, similar error bounds can
be derived for the $\ell_1$-regularized MLE and LSE, by following
almost verbatim the steps in \cite{chi:09}.  For brevity, we shall
omit the detail.  Instead, we shall focus on how to get error bounds
under the conditions in Proposition \ref{prop:basic2}.

\subsection{Exponential linear models}
Let $\{p(x;t): t\in I\}$ be a family of probability densities with
respect to a nonzero Borel measure $\mu$ on $\Reals$,
where  $I\subset\Reals$ is a closed interval, such that
$$
p(x;t)=\exp\Cbr{t y - \Lambda(t)},
\ \text{with} \ 
\Lambda(t) = \ln \Sbr{\int e^{t y}\,\mu(dy)}, \quad
t\in I.
$$

Suppose $\eno y n$ are independent, each with density $p(x;
X_i\tp\beta)$.  Let $D=\cD(I)$, where $\cD(I)$ is defined in
\eqref{eq:domain}.  Assume $\beta\in D$,  i.e. $X_i\tp \beta\in I$ for
each $i$.  The $\ell_1$-regularized MLE for $\beta$ is
\begin{align*}
  \est\beta = \argx\max_{v\in \cD(I)}
  \Sbr{
    y\tp Xv - \sum_{i=1}^n \Lambda(X_i\tp v) - c_r\normx{v}_1
  }.
\end{align*}

Let $\r_i = y_i - \mean(y_i) = y_i - \Lambda'(X_i\tp\beta)$,
$G(x)=\sum_{i=1}^n x_i$, $\psi_i(z) = \Lambda(z) -
\Lambda'(X_i\tp\beta) z$, and $\varphi_i(z) =z/2$.  Then it can be
been that $\est\beta$ satisfies the inequality \eqref{eq:ineq}.

Following almost verbatim the proof of Lemma 6.1 in \cite{chi:09}, 
if $\r$ satisfies the tail condition \eqref{eq:tail},
then Condition \ref{cond:A1} is satisfied by setting $c_0=c_\r$ and 
$$
c_1 = \sigma \sqrt{\frac{\ln(p/q)}{2n}}
\max_{\btw j 1 p}\normx{V_j}_2.
$$

On the other hand, in \cite{chi:09}, it was actually also shown that
for each $v\in\cD(I)$, $G(\psi(X v) - \psi(X\beta)) \ge (1/2)
\inf_{t\in I} \Lambda''(t) \times \normx{X(v-\beta)}_2^2$.  As a
result, we can set
$$
c_3 = (1/2) \inf_{t\in I} \Lambda''(t).
$$
If $\inf_{t\in I}\Lambda''(t)>0$, then, provided \eqref{eq:t} in
Proposition \ref{prop:basic2} is satisfied,
\begin{multline*}
  \prob\Biggl\{\normx{\est\beta - \beta}_2
  \le
  \frac{3(2+1/\tau) \sqrt{2+(1+2\tau)^2}}
  {a_X+b_X\mu_X} \times
  \sqrt{2\ln(p/q)}\\[-1ex]
  \times
  \frac{\sigma\sqrt{|\sppt(\beta)|}\max_{\btw j 1 p}\normx{V_j}_2}
  {n\inf_{t\in I}\Lambda''(t)}\Biggr\}\ge 1-c_\r q.
\end{multline*}
In particular, for the logistic model, where $\Lambda(t) = \ln
(1+e^t)$, since $\r_i = y_i - \Lambda'(X_i\tp\beta)$ with $y_i=0$ or
1, we can set $\sigma=1/2$ by Hoeffding's inequality
\cite{pollard:84}.  Furthermore, by $\Lambda''(t) =
(2\cosh(t/2))^{-2}$, $\inf_{t\in I} \Lambda''(t)>0$ for bounded $I$.

\subsection{Analytic models}
Suppose $y = f(X\tp\beta)+\r$, where $\r=(\eno\r n)\tp$ has mean 0 and
$f$ is defined on a closed interval $I\subset\Reals$ with positive
length.  Also, suppose $f$ can be continuously extended into an
analytic function on an open domain $\cN\subset\Coms$ that contains
$I$.   Now let $D\subseteq\cD(I)$ and assume $\beta\in D$.  The
$\ell_1$-regularized LSE estimator for $\beta$ is
\begin{align*}
  \est\beta = \argx\min_{v\in D}
  \Sbr{
    \normx{y-f(Xv)}_2^2 + c_r\normx{v}_1
  }.
\end{align*}
If we set $G(x) = \normx{x}_2^2$ and $\psi_i(z) = \varphi_i(z) =
f(z)$, then it can be seen that $\est\beta$ satisfies \eqref{eq:ineq},
and for $v\in I$, $G(\psi(Xv) - \psi(X\beta))\ge \sfd(f,I)^2
\normx{X(v-\beta)}_2^2$ \cite{chi:09}, where
\begin{align*}
  \sfd(f,I)=
  \inf \Cbr{\frac{|f(x)-f(y)|}{|x-y|}:\, x\in I,\, y\in I,\,
    x\not=y}.
\end{align*}
Therefore, if $\sfd(f,I)>0$, then we can set $c_3 = \sfd(f, I)^2$.

In order to apply Proposition \ref{prop:basic2}, we also need to get
$c_1$ for Condition \ref{cond:A1}.  We consider two cases.

In the first case, $D=\cD(I)\cap \{v\in \parsp:
\normx{v}_{1,\infty}\le\theta\rc/2\}$ and is compact, where $\theta\in
(0,1)$, $\rc>0$ such that $\{z\in \Coms: |z|\le \rc\}\subset\cN$, and 
$$
\normx{v}_{1,\infty} = \sum_{j=1}^p |v_j|\normx{V_j}_\infty.
$$
Let $\sigma$ be as in the tail condition \eqref{eq:tail}.  Given
$q\in (0,1)$, let $\lambda_p = \ln[p(1+q^{-1})]$.  As stated in
Proposition 6.5 in \cite{chi:09}, we can set
\begin{align*}
  c_1 = \sigma\sqrt{2\lambda_p}
  \sum_{k=1}^\infty
  \Sbr{
    \frac{\sqrt{k}|f\Sp k(0)|}{(k-1)!} (\theta\rc)^{k-1}\times
    n^{-\nth{2k}}\max_{\btw j 1 p} \normx{V_j}_{2k}
  }.
\end{align*}
Then by Proposition \ref{prop:basic2}, we get an error bound of the
same order as the $\ell_0$-regularized estimator in \cite{chi:09}.
Note that the constraints on $D$ include a bound on the weighted
$\ell_1$-norm $\normx{v}_{1,\infty}$ but no limits on $|\sppt(v)|$.
As a result, the LSE is purely regularized by $\ell_1$-norms
$\normx{v}_1$ and $\normx{v}_{1,\infty}$.

Second, $D=\cD(I)$ and is compact, but not
necessarily contained in a disc on which $f$ is analytic. 
Again, the LSE is purely regularized by $\ell_1$-norms of $v$.
However, it becomes harder to set $c_1$.  A relatively
simple choice of $c_1$ is as follows.  Let $\rc>0$, such that 
for any $x\in I$, $\{z\in \Coms: |z-x|<\rc\}\subset \cN$.  Let $d_k
= \sup_{x\in I} |f\Sp k(x)|/k!$, and $\delta(D)$ be the infimum of the
radii of spheres under $\normx{\cdot}_{1,\infty}$ that contain $D$,
i.e.,
$$
\delta(D) = \inf\{a>0: \text{there is $u\in\parsp$ such that }
\normx{v-u}_{1,\infty}<a \text{ for all } v\in D\}.
$$
Then, given $\rc_1\in (0, \rc)$, we can set
\begin{align} \label{eq:analytic}
  c_1 = \sqrt{2}\sigma \sum_{k=1}^\infty 
  \Sbr{
    k\sqrt{2p\ln(p Q)+k\lambda_p}\,d_k\rc_1^{k-1}
    \times n^{-\nth{2k}} \max_{\btw j 1 p} \normx{V_j}_{2k}
  },
\end{align}
where $Q=4\diam(D)/\rc_1+1$.  This value of $c_1$ results from
Proposition 5.5 (2) in \cite{chi:09} by noting the trivial bound
$|\sppt(v)|\le p$, which is nevertheless the tightest we can get, as 
no explicit constraints on $|\sppt(v)|$ are available.

Unfortunately, if we use \eqref{eq:analytic} to set $c_1$, then, in
order for the error bound in Proposition \ref{prop:basic2} to be at
most of order $o(1)$, $p$ cannot be very large.  Indeed, as the error
bound is proportional to $c_1\sqrt{|\sppt(\beta)|/n}\ge c
\sqrt{|\sppt(\beta)| p\ln p/n}$ for some $c>0$, $p$ has to be of order
$o(n/\ln n)$.

\subsection{Regression with noise-corrupted underlying linear
  structure}
It is possible to generalize the treatment for analytic models to
the following one
\begin{align} \label{eq:f-noise}
  y = f(X_i\tp\beta + \xi)+ \r
\end{align}
where $\eno \xi n$, $\eno\r n$ are independent with mean 0, and
$\xi_i$'s are identically distributed.  The model reflects the point
of view that noise can appear anywhere.  For nonlinear $f$, in general,
if the common distribution of $\xi_i$'s is unknown, then $\mean(y_i)$
are unknown and regression becomes impossible.  If, on the other hand,
the distribution is known, then $\mean[f(z+\xi_i)]$ are known.
Apprently, they are identical.  Denote $g(z)=\mean[f(z+\xi_1)]$ and
let $\delta_i = f(X_i\tp\beta+\xi_i) - g(X_i\tp\beta) + \r_i$.  Then 
\begin{align}  \label{eq:g}
  y = g(X\tp\beta)+ \delta.
\end{align}

Note that, in general, the distributions of $\delta_i$ depend
$X_i\beta$.  Since the latter are not identical, $\eno \delta n$ are
not identically distributed.  Furthermore, since $\beta$ is unknown,
in general, even if the distributions of $\r_i$ are known, the
distributions of $\delta_i$ are still unknown.  Despite this, by only
using the fact that $\delta_i$ are independent, each with mean 0, it
is possible to apply the results in previous sections to \eqref{eq:g},
hence getting error bounds of estimation for \eqref{eq:f-noise}.

To make this work, we need to check a few conditions, such as the
analyticity of $g(z)$ and the tail condition \eqref{eq:tail} for 
$\delta$.  We next present a case where the necessary conditions are
satisfied.

Suppose we set $D=\cD(I)$ with $I=[-R, R]$.  Suppose
$\xi_i$ are bounded random variables with $|\xi_i|<r$ and there is
$R_0>R+r$, such that $f$ is continuous on $\Delta_0:=\{z\in \Coms:
|z|\le R_0\}$ and analytic within it.  Let $\Delta=\{z\in\Coms:
|z|<R_0-r\}$.  For each $z\in\Delta$, by $z+\xi_1\in\Delta_0$,
$|f(z+\xi_1)|\le \sup_{\Delta_0} |f|<\infty$, so $g(z) =
\mean[f(z+\xi_1)]$ is well defined.  Clearly, $I$ is contained within
$\Delta$.

\begin{prop} \label{prop:g}
  (1) $g(z)$ is analytic on $\Delta$ and $\sfd(g,I)\ge
  \sfd(f,[-R_0,R_0])$.

  (2) If $\eno\r n$ satisfy \eqref{eq:tail} for some $\sigma>0$ and
  $c_\r>0$, then $\eno\delta n$ satisfy \eqref{eq:tail} as well for
  possibly different values of $\sigma$ and $c_\r$. 
  Moreover, if $\r_i$ are bounded, then $c_\r$ can always be
  set at 2.
\end{prop}

Thus, the results on $\ell_1$-regularized LSE in previous sections can
be applied to \eqref{eq:g}.  We omit the detail and will only prove
the Proposition.

\begin{proof}
  (1) Given $z\in\Delta$, for every possible value of $\xi_1$, we have
  $f(z+\xi_1) = \sum_{k=0}^\infty f\Sp k(\xi_1)z^k/k! $.  By Cauchy's
  contour integral,
  $$
  \frac{|f\Sp k(\xi_1)|}{k!}
  \le \nth{2\pi} \int_{|\zeta|=R_0} \frac{|f(\zeta)|
    d\zeta}{|\zeta-\xi_1|^{k+1}} 
  \le \frac{R_0 \sup_{\Delta_0} |f|}{(R_0-r)^{k+1}}
  $$
  Because $R_0-r>|z|$,
  $$
  \sum_{k=0}^\infty \frac{\mean |f\Sp k(\xi_1)|}{k!} |z|^k
  \le
  \frac{R_0\sup_{\Delta_0}|f|}{R_0-r}
  \sum_{k=0}^\infty \Grp{\frac{|z|}{R_0-r}}^k < \infty.
  $$
  Then by dominated convergence, it is seen that $g(z) =
  \sum_{k=0}^\infty \mean [f\Sp k(\xi_1)]z^k/k!$, with the power
  series being convergent on $\Delta$.   Therefore $g(z)$ is analytic
  on $\Delta$.
  
  To get $\sfd(g, I)\ge \sfd(f, [-R_0, R_0])$, let the right hand side
  be positive.  Then $f$ is monotone on $[-R_0, R_0]$, say,
  increasing.  Then $g(z)=\mean[f(z+\xi_1)]$ is increasing on $I$ and
  for $x<y$, $g(y)-g(x) = \mean[f(y+\xi_1) - f(x+\xi_1)]\ge \sfd(f,
  [-R_0, R_0])(y-x)$, finishing the proof of (1).

  (2)  Let $\eta_i=f(X_i\tp\beta+\xi_i)- g(X_i\tp\beta)$.  Then
  $\mathop{\rm ess\,sup} \eta_i - \mathop{\rm ess\,inf} \eta_i \le
  2\sup_{\Delta_0}|f|$ and $\delta_i = \eta_i+\r_i$.  Given
  $t\ge 0$ and $a\in\obssp$,
  \begin{align*}
    \prob\{|a\tp \delta| > t \normx{a}_2\}
    &
    \le \prob\{|a\tp \eta|> t\normx{a}_2/2\}
    +
    \prob\{|a\tp \r|> t\normx{a}_2/2\}
    \\
    &
    \le
    2 \exp\Cbr{
      -\frac{t^2}{8\sup_{\Delta_0}|f|^2}
    } + c_\r \exp\Cbr{
      -\frac{t^2}{8\sigma^2}
    },
  \end{align*}
  where the last inequality is due to Hoeffding's inequality and the
  tail condition \eqref{eq:tail}.  This implies the first claim of
  (2).  If $\r_i$ are bounded, then $\delta_i$ are bounded, and the
  second claim follows from Hoeffding's inequality.
\end{proof}

\bibliographystyle{acmtrans-ims2}
\bibliography{ldpdb,Estimate,NSdb.bib}

\end{document}

%% file: myamsdefs
\long\def\comment#1{}

\def\qed{\hspace*{\fill}$\Box$\medbreak}

\def\vf#1{\boldsymbol{#1}}

\def\Grp#1{\left(#1\right)}

\def\Cbr#1{\left\{#1\right\}}

\def\Sbr#1{\left[#1\right]}

\def\Abs#1{\left|#1\right|}

\def\Ang#1{\left\langle#1\right\rangle}

\def\cf#1{\mathbf{1}\Cbr{#1}}

\def\eno#1#2{{#1}_1, \ldots, {#1}_{#2}}         

\def\nth#1{\frac{1}{#1}}

\def\Sp#1{\sp{(#1)}}

\def\Coms{\mathbb{C}}

\def\Reals{\mathbb{R}}

